\documentclass[a4paper,fleqn,10pt,twocolumn]{SICE_ISCS}
\usepackage{amsmath,amssymb,latexsym,mathtools}
\usepackage[utf8]{inputenc}
\usepackage{graphicx}
\usepackage[caption=false,font=footnotesize]{subfig}\usepackage{xcolor,soul}
\usepackage[colorlinks,urlcolor=blue,linkcolor=blue,citecolor=blue]{hyperref}
\usepackage{algorithm2e}
\usepackage{natbib}
\usepackage{jabbrv}

\setcitestyle{square,numbers}

\newcommand{\htwo}{\ensuremath{\text{H}_2} }
\newcommand{\heh}{\ensuremath{\text{HeH}^+} }
\newcommand{\lih}{\ensuremath{\text{LiH}} }

\DeclareMathOperator{\tr}{trace}
\DeclareMathOperator{\vvec}{vec}
\DeclareMathOperator{\unvec}{unvec}
\DeclareMathOperator{\diag}{diag}

\title{Nonlinear Optimal Control of Electron Dynamics within Hartree-Fock Theory}

\author{Harish~S. Bhat${}^{1,2\dagger}$ and Hardeep Bassi${}^{1}$ and Christine~M. Isborn${}^{3}$}
\speaker{Harish~S. Bhat}

\affils{$^{1}$Department of Applied Mathematics, University of California, Merced, Merced, CA 95343 USA\\
${}^{2}$E-mail: hbhat@ucmerced.edu\\
$^{3}$Department of Chemistry and Biochemistry, University of California, Merced, Merced, CA 95343 USA\\
}

\abstract{%
Consider the problem of determining the optimal applied electric field to drive a molecule from an initial state to a desired target state.  For even moderately sized molecules, solving this problem directly using the exact equations of motion---the time-dependent Schr\"odinger equation (TDSE)---is numerically intractable.  We present a solution of this problem within time-dependent Hartree-Fock (TDHF) theory, a mean field approximation of the TDSE.   Optimality is defined in terms of minimizing the total control effort while maximizing the overlap between  desired and achieved target states.  We frame this problem as an optimization problem constrained by the nonlinear TDHF equations; we solve it using trust region optimization with gradients computed via a custom-built adjoint state method.  For three molecular systems, we show that with very small neural network parametrizations of the control, our method yields solutions that achieve desired targets within acceptable constraints and tolerances.
}


\begin{document}
\maketitle

\section{Introduction}
The ability to control electron dynamics enables control of chemical reactions \cite{Assion1998,Rabitz2001}, energy flow in light harvesting complexes \cite{Herek2002}, the aromaticity of benzene \cite{ulusoyCorrelatedElectronDynamics2011}, and vacancy centers in diamond \cite{Tian2019}, among other examples \cite{Werschnik_2007,Keefer2018}.  For even moderately sized molecular systems, an obstacle to solving such control problems is our inability to solve the time-dependent Schr\"odinger equation (TDSE) without the use of approximations and numerical methods.  Time-dependent Hartree-Fock (TDHF) theory, with the governing equation
\begin{equation}
\label{eqn:tdhf}
i \frac{dP}{dt} = [H(P(t),t),P(t)],
\end{equation}
is one computationally tractable approximation of the electronic TDSE.  Here $[H,P] = HP - PH$ denotes the matrix commutator, $H$ denotes the Hamiltonian (defined in more detail below), and $P(t)$ is the time-dependent $1$-electron reduced density matrix from which many observables of interest can be computed.  Let $P_0$ and $P_T$ denote desired initial and final states.  \emph{What is the optimal electric field  that we must apply to a molecular system governed by (\ref{eqn:tdhf}) to take it from initial state $P_0$ to final state $P_T$?} In this paper, we propose a solution to this that blends classical model-based optimal control with a neural network parametrization of the controller.

The derivation of (\ref{eqn:tdhf}) from the TDSE is standard in the literature \cite{mcweeny1989methods,szabo2012modern}; for a conceptual overview, see Section \ref{sect:background}. In this derivation, the TDSE's Hamiltonian operator $\hat{H}$ becomes the density-dependent Hartree-Fock Hamiltonian matrix $H$ in (\ref{eqn:tdhf}).  \emph{Because $H$ depends on $P$, the right-hand side of (\ref{eqn:tdhf}) is nonlinear in $P$}.  In our review of the literature below, we find that optimal control problems for (\ref{eqn:tdhf}) have not been considered in prior work.

The $N \times N$ Hermitian matrices $H$ and $P$ satisfy
\begin{equation}
\label{eqn:hfham}
H(P,t) = H_0(P) + V^{\text{ext}}(t),
\end{equation}
where $H_0(P)$ is the field-free Hamiltonian and $V^{\text{ext}}(t)$ is an external, time-dependent potential, e.g., an applied electric field that can drive the system between states.  The field-free Hamiltonian $H_0(P)$ includes the kinetic energy of the electrons, the potential due to interaction with the fixed nuclei, and the density-dependent inter-electronic repulsion term, also known as the Fock matrix in Hartree-Fock theory.

We work in the dipole approximation: when we apply an electric field with real amplitudes ${\mathbf{a}}(t) = [a_1(t),a_2(t),a_3(t)]$ in the $x$, $y$, and $z$ directions, the corresponding potential is
\begin{equation}
\label{eqn:actuation}
V^{\text{ext}}(t) = \sum_{j=1}^{3} a_j(t) M_j.
\end{equation}
Here $\{M_1, M_2, M_3\}$ are the $N \times N$ dipole moment matrices in the $x$, $y$, and $z$ directions, respectively.  These Hermitian matrices depend on the molecular system at hand.  One can see immediately that as $N$ increases, the system becomes severely underactuated.  Let ${\mathcal{S}}(N)$ denote the vector space of $N \times N$ complex Hermitian matrices; this space has real dimension $N^2$.  We see that (\ref{eqn:actuation}) constrains $V^{\text{ext}}(t)$ to lie in at most a three-dimensional subspace of ${\mathcal{S}}(N)$.  For 2-electron linear diatomic molecular systems that we study in this paper, we will find that $M_1 = M_2 \equiv 0$, constraining $V^{\text{ext}}(t)$ still further.

Purely as a mathematical/computational device, we consider feedback control whereby the amplitude ${\mathbf{a}}(t)$ in (\ref{eqn:actuation}) is a function of the electron density at time $t$:
\begin{equation}
\label{eqn:feedback}
{\mathbf{a}}(t) = {\mathbf{a}}(P(t)).
\end{equation}
By device, what is meant is that after solving for the feedback control ${\mathbf{a}}(P)$, we use it to simulate (\ref{eqn:tdhf}) with (\ref{eqn:hfham}), (\ref{eqn:actuation}), and (\ref{eqn:feedback}).  While doing so, we record the left-hand side of (\ref{eqn:feedback}), the open-loop control ${\mathbf{a}}(t)$ that would have produced the desired trajectory.  \emph{The goal of this work is to demonstrate that very small neural networks can produce reasonable control solutions at  low computational cost.}


\subsection{Prior Work}
The equation of motion (\ref{eqn:tdhf}) differs substantially from prior Liouville-von Neumann (LvN) equations for which optimal control problems have been considered previously \cite{Borzi2017,goodwinAdvancedOptimalControl2017,dehaghaniQuantumOptimalControl2022,dehaghaniHighFidelityQuantum2022,Aguiar2023,goodwinAcceleratedNewtonRaphsonGRAPE2023}.  For any quantum system, suppose we form the density operator and express it in an orthonormal basis as a matrix $P(t)$.  Then starting from the TDSE, we can easily derive the LvN equation
\begin{equation}
\label{eqn:lvn}
i \frac{dP}{dt} = [H(t), P(t)].
\end{equation}
This differs from (\ref{eqn:tdhf}) in two key aspects.  In (\ref{eqn:tdhf}), $P(t)$ refers to the \emph{1-electron reduced density matrix} and the Hamiltonian is density-dependent; in (\ref{eqn:lvn}), $P(t)$ is the full density matrix and the Hamiltonian does not depend on density.  In short, the dynamics of (\ref{eqn:tdhf}) are nonlinear even when we consider no applied field or control, a regime in which (\ref{eqn:lvn}) has linear dynamics.

One can formulate TDHF theory in real three-dimensional space, rather than in a finite basis set as in (\ref{eqn:tdhf}).  The result is a nonlinear, nonlocal partial differential equation in ${\mathbb{R}}^3$; in this framework, an optimal control problem for the Helium atom has been studied \cite{bonnard2012optimal}. The finite basis TDHF equation (\ref{eqn:tdhf}) is commonly used in the quantum/computational chemistry literature \cite{Li2005,Isborn2008,Isborn2009,Ranka2023}.

Time-dependent density functional theory (TDDFT), which also  features a density-dependent Hamiltonian, is another widely used and computationally tractable approximation of the TDSE.  Thus prior studies of  optimal control for TDDFT share structural similarities with the present work  \cite{castroControllingDynamicsManyElectron2012,sprengelInvestigationOptimalControl2018,Wong2024CPC}.  However, the equations of motion for TDDFT consist of nonlinear, nonlocal partial differential equations in ${\mathbb{R}}^3$; our approach is mathematically and computationally simpler.

Below we will make use of the adjoint method and alternating forward/backward propagation steps, which have been used to solve optimal control problems for the TDSE \cite{Rabitz1988,Goerz2019,RAZA2021107541,Wong2024Travolta,Shao2024}.  Without approximations of some kind (such as Hartree-Fock theory), such approaches cannot scale to the Hamiltonians that govern electron dynamics in molecules.

We parameterize the control signal using a neural network, and fuse this neural network model with the governing equations (\ref{eqn:tdhf}-\ref{eqn:feedback}).  This process leverages known physics, and thus differs from recent approaches (for the TDSE) that use neural networks to effectively replace the equations of motion \cite{Wong2020PCCP,dehaghaniApplicationPontryaginNeural2023}.


\subsection{Key Contributions}
In this paper, we give the first solution to an optimal control problem for TDHF dynamics governed by (\ref{eqn:tdhf}).  We discretize (\ref{eqn:tdhf}) with a symmetry-preserving numerical method and we parameterize the feedback control ${\mathbf{a}}(P)$ using a neural network.  Our method fuses (i) an adjoint method that differentiates through the symmetry-preserving discretization of (\ref{eqn:tdhf}) with (ii) backpropagation to compute required Jacobians of the neural network.  We provide empirical evidence that for three molecular systems, even with tiny neural networks, our methods yield low-cost control inputs ${\mathbf{a}}(t)$ that drive the dynamics to final states $P(T)$ that closely match desired targets $P_T$.  By using such small neural networks for these three model systems, we show that our methods have a chance to scale to much larger systems.
\section{Problem}

\subsection{Background and Continuous-Time Formulation}
\label{sect:background}
To derive (\ref{eqn:tdhf}) starting from the TDSE, we begin by holding fixed all particles other than electrons.  As part of the Born-Oppenheimer approximation, we consider only the time-evolution of the electronic wave function $\Psi$ in the field of fixed nuclei.  Next, we approximate $\Psi$ using a single Slater determinant computed from a finite set of spin orbitals  \cite{szabo2012modern}.  These spin orbitals are products of spin variables with spatial functions (molecular orbitals) that are themselves linear combinations of atomic orbitals (AOs).  The number $N$ of AOs (i.e., the size of the \emph{basis set})  is up to the user: for a given molecule, increasing $N$ leads to more accurate models.    Applying these approximations to the TDSE, and expressing all operators in an orthonormal basis, one can derive (\ref{eqn:tdhf}) as the governing equation for $P(t)$, the time-dependent $1$-electron reduced density matrix that corresponds to $\Psi$ \cite{mcweeny1989methods,szabo2012modern}. 

Suppose our system has $N_{\text{e}}$ total electrons and let ${\mathcal{S}}(N)$ denote the space of Hermitian $N \times N$ matrices.  One can contrast the TDSE, a linear partial differential equation in an intractably high-dimensional space (ignoring spin, the space will be ${\mathbb{R}}^{3 N_{\text{e}}}$), with the TDHF equation (\ref{eqn:tdhf}), a tractable, nonlinear system of ordinary differential equations in the space ${\mathcal{S}}(N)$.

Let ${\mathcal{A}} \subset \{1, 2, 3\}$ denote the set of active dipole moment matrices in (\ref{eqn:actuation}); let $N_{a} = |{\mathcal{A}}|$.  From (\ref{eqn:feedback}), our control input is ${\mathbf{a}} : {\mathcal{S}}(N) \to {\mathbb{R}}^{N_{a}}$; with this, (\ref{eqn:actuation}) is
\begin{equation}
\label{eqn:feedbackactuation}
V^{\text{ext}}(P(t)) = \sum_{j \in {\mathcal{A}}} a_j(P(t)) M_j.
\end{equation}
Because of the dimensionality considerations mentioned above, and also because one or more of the $M_j$ matrices may be rank-deficient, the system is underactuated.  As (\ref{eqn:tdhf}) is nonlinear, checking its controllability is nontrivial; in numerical experiments, we find that our ability to steer the system to the terminal state $P_T$ depends both on the value of $P_T$ and the time $T > 0$ that we choose.  Therefore, instead of imposing $P_T$ as an equality constraint ``$P(T) = P_T$,'' we instead  \emph{maximize the fidelity} ${\mathcal{F}}$ between the achieved state $P(T)$ and the target $P_T$:
\begin{equation}
\label{eqn:overlap}
{\mathcal{F}}(P(T), P_T) = \tr \bigl(P(T) P_T P(T)\bigr).
\end{equation}
This is the $p=2$ case of a physically motivated norm $\| \cdot \|_p$ that enjoys properties such as symmetry, nondegeneracy, and invariance under unitary transformations \cite{Liang2019}; the $p=1$ version of this norm has been used in prior optimal control studies \cite{dehaghaniHighFidelityQuantum2022,dehaghaniQuantumOptimalControl2022,Aguiar2023}.  To realize the equivalence with the definition in \cite{Liang2019}, one must recall that electron density matrices $P(t)$ that satisfy (\ref{eqn:tdhf}) are Hermitian, idempotent, and have constant trace equal to $N_{\text{e}}/2$, half the number of electrons in the system.  Using these properties, we can show that $0 \leq {\mathcal{F}}(P(T), P_T) \leq N_{\text{e}}/2$.

With (\ref{eqn:feedbackactuation}) and (\ref{eqn:overlap}), the continuous-time formulation of our optimal control problem is to find ${\mathbf{a}} : {\mathcal{S}}(N) \to {\mathbb{R}}^{N_{a}}$ that minimizes the objective
\begin{equation}
\label{eqn:objective}
\frac{1}{2} \int_{t=0}^{t=T} \left\| V^{\text{ext}}(P(t)) \right\|_F^2 \, dt - \frac{\rho}{2} {\mathcal{F}}(P(T), P_T)^2
\end{equation}
subject to the equation of motion (\ref{eqn:tdhf}) and $P(0) = P_0$.  Here $P_0$ and $P_T$ are prescribed initial and final states.  We have formed this objective by combining a running cost on the squared norm of the applied field together with a terminal cost, the lack of fidelity between the final achieved state and the target.  Since the running cost has no upper bound, the user-defined factor $\rho > 0$ exists to balance the running and terminal costs.



\subsection{Discrete-Time Dynamics} To solve the optimal control problem described above, we pursue a discretize-then-optimize strategy.  We begin by discretizing (\ref{eqn:tdhf}) via the widely-used modified midpoint unitary transform (MMUT) method \cite{Li2005}.  We fix the step size $\Delta t > 0$ and use $P^k$ to denote our numerical approximation of $P(k \Delta t)$.  In the MMUT scheme, at time step $k=0$, we compute $P^1$ from $P^0$:
\begin{equation}
\label{eqn:MMUTfirststep}
P^1 = U^0 P^0 {U^0}^\dagger \ \text{ with } \ U^0 = \exp(-i \Delta t H(P^0, 0)).
\end{equation}
At subsequent time steps $k = 1, 2, \ldots, K-1$, we compute $P^{k+1}$ using both $P^k$ and $P^{k-1}$:
\begin{subequations}
\label{eqn:MMUTothersteps}
\begin{align}
U^{k} &= \exp(-2i \Delta t H(P^k, k \Delta t)) \\
P^{k+1} &= U^k P^{k-1} {U^k}^\dagger.
\end{align}
\end{subequations}
This scheme preserves Hermitian symmetry, idempotency, and constant trace of $P^k$ for all $k$. 
\section{METHODS}
\label{sect:methods}

\subsection{Real Representation and Control}
\label{sect:controlmodel}
For $P \in {\mathcal{S}}(N)$, Hermitian symmetry means that the matrix $P$ is determined by (i) the $N(N+1)/2$ real parts of its upper-triangular entries (including the diagonal) and (ii) the $N(N-1)/2$ imaginary parts of its lower-triangular entries (excluding the diagonal).  In total, the matrix $P$ is determined by a  vector ${\mathbf{p}} \in {\mathbb{R}}^{N^2}$.  We now explain how to pass systematically between $P$ and ${\mathbf{p}}$.  First, using (i) and (ii), we can easily form a \emph{basis} for the space ${\mathcal{S}}(N)$.  When $N = 2$, the basis has $N^2 = 4$ elements:
\begin{equation}
\label{eqn:explicitbasis}
\begin{bmatrix} 1 & 0 \\ 0 & 0 \end{bmatrix},
\begin{bmatrix} 0 & 1 \\ 1 & 0 \end{bmatrix},
\begin{bmatrix} 0 & 0 \\ 0 & 1 \end{bmatrix},
\begin{bmatrix} 0 & i \\ -i & 0 \end{bmatrix}.
\end{equation}
Let $B_{q,r,m}$ denote this basis, so that $B_{:,:,m}$ is the $m$-th $N \times N$ matrix given above, for $m = 1, 2, \ldots, N^2$.  Given $P \in {\mathcal{S}}(N)$, we have
\begin{equation}
\label{eqn:basisrep}
P_{q,r} = \sum_{m=1}^{N^2} B_{q,r,m} p_m ,
\end{equation}
Given any complex $N \times N$ matrix $P$, let $\vvec(P) \in {\mathbb{C}}^{N^2}$ denote its vectorized representation, the concatenation of the rows of $P$: $\vvec(P)_{(i-1)N + j} = P_{i,j}$.  Note that $\vvec(P)$ collapses two indices into one.  Applying $\vvec$ to collapse the $(q,r)$ indices on both sides of (\ref{eqn:basisrep}), we obtain
$\vvec(P) = R {\mathbf{p}}$,
where $R$ is an $N^2 \times N^2$ matrix such that $R_{(q-1)N + r, m} = B_{q,r,m}$.  The matrix $R$ transforms the real representation ${\mathbf{p}}$ into the vectorized version of the matrix $P$.  If we reshape $\vvec(P)$ to size $N \times N$, we recover $P \in {\mathcal{S}}(N)$.  Note that $R$ is invertible and hence we can go back: ${\mathbf{p}} = \unvec(P) := R^{-1} \vvec(P)$.

With the above, we can replace ${\mathbf{a}}(P)$ with ${\mathbf{a}}({\mathbf{p}})$.  We parameterize this using a neural network with parameters $\theta \in {\mathbb{R}}^{M}$ and write ${\mathbf{a}}({\mathbf{p}}; \theta)$.  In effect, we have replaced ${\mathbf{a}}(P)$, a function of a complex Hermitian matrix, with ${\mathbf{a}} : {\mathbb{R}}^{N^2} \times {\mathbb{R}}^{M} \to {\mathbb{R}}^{N_a}$,
a function with real inputs and real outputs.  We assume access to Jacobians of ${\mathbf{a}}$ with respect to either ${\mathbf{p}}$ or $\theta$, computed via automatic differentiation.  When we write $V^{\text{ext}}(P; \theta)$ in what follows, we mean (\ref{eqn:actuation}) in conjunction with ${\mathbf{a}}(t) = {\mathbf{a}}(\unvec{P(t)}; \theta)$.

\subsection{Discrete-Time Optimal Control Problem}
\label{sect:lagrangian}
With the discretized scheme and the neural network control function, the \emph{discrete-time optimal control problem} is to find $\theta$ that minimizes the objective
\begin{equation}
\label{eqn:discobj}
{\mathcal{J}}(\theta) = \frac{1}{2} \sum_{k=0}^{K-1} \left\| V^{\text{ext}}(P^k; \theta) \right\|_F^2 - \frac{\rho}{2} {\mathcal{F}}(P^K, P_T)^2
\end{equation}
subject to the dynamics (\ref{eqn:MMUTfirststep}-\ref{eqn:MMUTothersteps}) and $P^0 = P_0$.  We introduce the Frobenius inner product on the space of complex $n \times n$ matrices:
\begin{equation}
\label{eqn:froip}
\langle A, B \rangle = \tr A^\dagger B.
\end{equation}
Properties of this inner product are detailed in the Appendix.  With this, we form the real-valued Lagrangian
\begin{multline}
\label{eqn:Lagrangian}
L(P, \Lambda, \theta) = \frac{1}{2} \sum_{k=0}^{K-1} \left\| V^\text{ext}(P^k; \theta) \right\|_F^2  - \frac{\rho}{2} {\mathcal{F}}(P^K, P_T)^2 \\
- \Re \langle \lambda^1, P^1 - U^0 P_0 {U^0}^\dagger \rangle   \\
- \Re \sum_{k=1}^{K-1} \langle \lambda^{k+1}, P^{k+1} - U^k P^{k-1} {U^k}^\dagger \rangle.
\end{multline}
In this Lagrangian and in what follows, $\{U^k\}$ refers to the unitary matrices defined in (\ref{eqn:MMUTfirststep}) and (\ref{eqn:MMUTothersteps}); also, we identify $P^0$ with $P_0$. We formed (\ref{eqn:Lagrangian}) by combining the objective (\ref{eqn:discobj}) with terms that enforce the dynamics as constraints.  We refer to the Lagrange multipliers $\Lambda = \{\lambda^k\}_{k=1}^{K}$ as adjoint variables.

\subsection{Discrete-Time Adjoint Equation}
By the Lagrange multiplier theorem \cite{luenberger1997optimization}, if $\theta$ yields a feasible solution of the discrete-time optimal control problem, then there exists a corresponding critical point $(P, \Lambda, \theta)$ of the Lagrangian.  Hence we compute gradients.  As $\nabla_{\Lambda} L = 0$ gives us (\ref{eqn:MMUTfirststep}-\ref{eqn:MMUTothersteps}), we consider the gradient of $L$ with respect to $P$.  Using $:$ to denote tensor contraction as in (\ref{eqn:contractions}), we begin with
\begin{align*}
&\frac{1}{2} \frac{d}{d\epsilon}\biggr|_{\epsilon=0} \left\| V^\text{ext}(P^k + \epsilon \delta P^k; \theta) \right\|_F^2 \\
 &= \frac{1}{2} \bigg\langle \frac{dV}{dP}(P^k; \theta) : \delta P^k, V^\text{ext}(P^k; \theta) \bigg\rangle \\
 &\qquad\qquad\qquad+ \frac{1}{2} \bigg\langle V^\text{ext}(P^k; \theta), \frac{dV}{dP}(P^k; \theta) : \delta P^k \bigg\rangle \\
 &= \Re  \bigg\langle V^\text{ext}(P^k; \theta), \frac{dV}{dP}(P^k; \theta) : \delta P^k \bigg\rangle.
\end{align*}
Since the initial condition is fixed, $\delta P^0 \equiv 0$, and so
\begin{multline*}
\frac{d}{d\epsilon} L(P + \epsilon \delta P, \Lambda, \theta) \biggr|_{\epsilon=0} = - \Re \sum_{k=0}^{K-1} \langle \lambda^{k+1}, \delta P^{k+1} \rangle \\
+ \Re \sum_{k=1}^{K-1}  \biggl\langle V^\text{ext}(P^k; \theta), \frac{\partial V^\text{ext}}{\partial P} (P^k; \theta) : \delta P^k \biggr\rangle  \\
+ \Re \sum_{k=1}^{K-1} \langle \lambda^{k+1}, (\delta U^k) P^{k-1} {U^k}^\dagger \\
+ U^k (\delta P^{k-1}) {U^k}^\dagger + U^k P^{k-1} (\delta {U^k})^\dagger \rangle.
\end{multline*}
Here $\displaystyle \delta U^k = \frac{d}{d\epsilon}\biggr|_{\epsilon=0} U^k[P^k + \epsilon \delta P^k]$.
Applying (\ref{eqn:tens}) to the first term on the first line, reindexing the first term on the next line, and then using (\ref{eqn:mat12}),
\begin{multline}
\label{eqn:thread2}
\frac{d}{d\epsilon} L(P + \epsilon \delta P, \Lambda, \theta) \biggr|_{\epsilon=0} = - \Re \sum_{k=1}^{K} \langle \lambda^{k}, \delta P^{k} \rangle \\
+ \Re \sum_{k=1}^{K-1}  \biggl\langle V^\text{ext}(P^k; \theta) : \overline{\frac{\partial V^\text{ext}}{\partial P} (P^k; \theta)} , \delta P^k \biggr\rangle \\
+ \Re \sum_{k=1}^{K-1} \langle \lambda^{k+1} U^k P^{k-1}, \delta U^k \rangle \\
+ \Re \sum_{k=1}^{K-2} \langle (U^{k+1})^\dagger \lambda^{k+2} U^{k+1}, \delta P^{k} \rangle \\
+ \Re \sum_{k=1}^{K-1} \langle P^{k-1} (U^k)^\dagger \lambda^{k+1},  (\delta {U^k})^\dagger \rangle.
\end{multline}
To compute $\delta U^k$ and $(\delta U^k)^\dagger$, note that for $k \geq 1$,
\begin{equation}
\label{eqn:deltaUkDEF}
\delta U^k = \frac{d}{d\epsilon}\biggr|_{\epsilon=0} \! \! \! \! \! \exp(-2i \Delta t H(P^k + \epsilon \delta P^k, k \Delta t; \theta)).
\end{equation}
Via the chain rule,
\begin{multline*}
(\delta U^k)_{a,b} = \sum_{j,\ell} \frac{ \partial }{\partial Z_{j,\ell}} (\exp Z)_{a,b}  \biggr|_{Z = - 2 i \Delta t H(P^k, k \Delta t; \theta)} \\
\cdot (-2i \Delta t) \frac{d}{d\epsilon}\biggr|_{\epsilon=0} H_{j,\ell}(P^k + \epsilon \delta P^k, k \Delta t; \theta).
\end{multline*}
When restricted to Hermitian or anti-Hermitian matrices $Z$,  the Jacobian of the matrix exponential, $\partial (\exp Z)_{a,b} / \partial Z_{j,\ell}$, follows from prior results \cite{Lewis2001}.  In our code, we use an efficient implementation that stems from our own derivations; we omit these for reasons of space.  Continuing the derivation from above,
\begin{multline*}
\frac{d}{d\epsilon}\biggr|_{\epsilon=0} H_{j,\ell}(P^k + \epsilon \delta P^k, k \Delta t; \theta) \\
= \sum_{r,s} \frac{\partial V^\text{ext}_{j,\ell}}{\partial P_{r,s}}(P^k,k\Delta t; \theta) (\delta P^k)_{r,s},
\end{multline*}
which implies
\begin{gather}
\label{eqn:deltaUk}
\delta U^k = (-2i \Delta t) \zeta^{k,-} : \delta P^k \\
\zeta^{k,-}_{a,b,r,s} = \sum_{j,\ell} \frac{\partial (\exp Z)_{a,b}}{\partial Z_{j,\ell}} \Biggr|_{\substack{\\ \\ Z = - 2 i \Delta t H(P^k, k \Delta t; \theta)}} \! \! \! \! \! \! \! \! \! \! \! \! \! \! \! \! \! \! \! \! \! \! \! \! \!  \! \! \! \! \! \! \! \!  \cdot \frac{\partial V^\text{ext}_{j,\ell}}{\partial P_{r,s}}(P^k,k\Delta t; \theta) \nonumber,
\end{gather}
a quantity that we can evaluate using automatic differentiation (to find the Jacobian of ${\mathbf{a}}({\mathbf{p}})$) and linear algebra (to convert this to the Jacobian of $V^\text{ext}(P)$).  The complex conjugate can be dealt with similarly.  Note that $[\exp(Z)]^\dagger = \exp(Z^\dagger)$.  For $Z = -2i \Delta t H$, where $H$ is Hermitian, we have $Z^\dagger = 2i \Delta t H$.  Therefore, for $k \geq 1$, the definition of the derivative gives
\begin{equation}
\label{eqn:deltaUkdaggerDEF}
\delta {U^k}^\dagger 
= \frac{d}{d\epsilon}\biggr|_{\epsilon=0} \! \! \! \! \! \exp(2i \Delta t H(P^k + \epsilon \delta P^k, k \Delta t; \theta)).
\end{equation}
This means that changing the sign of $i$ in (\ref{eqn:deltaUk}) gives the correct result.  For $k \geq 1$, 
\begin{gather}
\label{eqn:deltaUkdagger}
\delta {U^k}^\dagger = (2i \Delta t) \zeta^{k,+} : \delta P^k \\
\nonumber \zeta^{k,+}_{a,b,r,s} = \sum_{j,\ell} \frac{\partial (\exp Z)_{a,b}}{\partial Z_{j,\ell}} \Biggr|_{\substack{\\ \\ Z = 2 i \Delta t H(P^k, k \Delta t; \theta)}} \! \! \! \! \! \! \! \! \! \! \! \! \! \! \! \! \! \! \! \! \! \! \! \! \!  \! \! \! \! \! \! \! \! \cdot \frac{\partial V^\text{ext}_{j,\ell}}{\partial P_{r,s}}(P^k,k\Delta t; \theta).
\end{gather}
There is one last derivative to compute, that of the fidelity with respect to the final state $P^K$.  Noting that the trace of a matrix is the sum of its eigenvalues, and using $\mu_j(Z)$ to denote the $j$-th eigenvalue of an $N \times N$ matrix $Z$,
\begin{align}
&\frac{\partial}{\partial P^K} \left( \frac{\rho}{2} {\mathcal{F}}(P^K,P_T)^2 \right) \nonumber \\
 &= \rho {\mathcal{F}}(P^K,P_T)  \sum_{j=1}^{N} \frac{ \partial \mu_j( Z ) } {\partial Z}  \Biggr|_{\substack{\\ \\ Z = P^K P_T P^K}} \! \! \! \! \! \! \! \! \! \! \! \! \! \! \! \! \! \! \! \! \! \! \! \! \! \! \! \cdot  \frac{\partial  (P^K P_T P^K)}{\partial P^K} \nonumber \\
\label{eqn:fidelityderiv}
 &= \rho {\mathcal{F}}(P^K,P_T)  \sum_{j=1}^{N} v_j v_j^\dagger P^K P_T + P_T P^K  v_j v_j^\dagger,
\end{align}
where $v_j$ denotes the $j$-th eigenvector of $P^K P_T P^K$. We now return to (\ref{eqn:thread2}).  Using (\ref{eqn:deltaUk}), (\ref{eqn:deltaUkdagger}), and (\ref{eqn:tens}), we obtain
\begin{multline}
\label{eqn:thread5}
\frac{d}{d\epsilon} L(P + \epsilon \delta P, \Lambda, \theta) \biggr|_{\epsilon=0} \! \! \! \! \! = -\Re \langle \text{result of (\ref{eqn:fidelityderiv})}+ \lambda^{K}, \delta P^K \rangle\\
+ \sum_{k=1}^{K-1}  \Re  
 \biggl\langle - \lambda^k + V^\text{ext}(P^k; \theta) : \overline{\frac{\partial V^\text{ext}}{\partial P} (P^k; \theta)} \\
 +  (U^{k+1})^\dagger \lambda^{k+2} U^{k+1} +  (2i \Delta t) \lambda^{k+1} U^k P^{k-1} : \overline{\zeta^{k,-}} \\
+ (-2i \Delta t) P^{k-1} (U^k)^\dagger \lambda^{k+1} :  \overline{\zeta^{k,+}}, \delta P^k \biggr\rangle,
\end{multline}
where as shorthand we have introduced $\lambda^{K+1} \equiv 0$.
For each fixed $k$ from $1$ to $K$, we want the right-hand side to vanish for \emph{all variations} $\delta P^k$.  With the Frobenius norm (\ref{eqn:froip}), one can show that if $X$ is an $n \times n$ complex matrix such that $\Re \langle X, Y \rangle = 0$ for all complex $n \times n$ matrices $Y$, then $X = 0$.  Thus we obtain the final condition
\begin{equation}
\label{eqn:fincond}
\lambda^K \! = \! -\rho {\mathcal{F}}(P^K \! ,  P_T)  \sum_{j=1}^{N} v_j v_j^\dagger P^K P_T + P_T P^K  v_j v_j^\dagger ,
\end{equation}
the first backward step
\begin{multline}
\label{eqn:firstlambdastep}
\lambda^{K-1} = V^\text{ext}(P^{K-1}; \theta) : \overline{\frac{\partial V^\text{ext}}{\partial P} (P^{K-1}; \theta)}  \\
+  2i \Delta t \lambda^{K} U^{K-1} P^{{K-2}} : \overline{\zeta^{K-1,-}} \\
- 2i \Delta t P^{K-2} (U^{K-1})^\dagger \lambda^{K} :  \overline{\zeta^{K-1,+}},
\end{multline}
and the backward-in-time adjoint equation
\begin{multline}
\label{eqn:nextlambdastep}
\lambda^k = V^\text{ext}(P^k; \theta) : \overline{\frac{\partial V^\text{ext}}{\partial P} (P^k; \theta)} +  (U^{k+1})^\dagger \lambda^{k+2} U^{k+1} \\
+  2i \Delta t \lambda^{k+1} U^k P^{k-1} : \overline{\zeta^{k,-}} \\
- 2i \Delta t P^{k-1} (U^k)^\dagger \lambda^{k+1} :  \overline{\zeta^{k,+}}
\end{multline}
for $k = K-2, \ldots, 1$.  Note that (\ref{eqn:fincond}-\ref{eqn:nextlambdastep}) can also be derived via Pontryagin's minimum principle \cite{d2007introduction,sontag2013mathematical}.
\subsection{Optimization Algorithm}
\label{sect:optalg}
Using the above results, we now explain how we solve the discrete-time optimal control problem stated at the top of Section \ref{sect:lagrangian}.  Suppose we have a current iterate $\theta^{(m)}$, the parameters of our neural network control.  Given $\theta^{(m)}$ and the initial condition $P^0 = P_0$, we solve (\ref{eqn:MMUTfirststep}) and (\ref{eqn:MMUTothersteps}) to generate a trajectory $P = \{P^k\}_{k=0}^K$.  Using this trajectory and $\theta^{(m)}$, we then use (\ref{eqn:fincond}), (\ref{eqn:firstlambdastep}) and (\ref{eqn:nextlambdastep}) to generate the adjoint trajectory $\Lambda = \{\lambda^k\}_{k=1}^{K}$.  By carrying out the steps until this point, we have found $(P,\Lambda,\theta)$ such that $\nabla_{P} L = 0$ and $\nabla_{\Lambda} L = 0$.  There is only one remaining gradient, which we now compute:
\begin{multline}
\label{eqn:dLdtheta}
\nabla_{\theta} L = \Re \sum_{k=0}^{K-1} \biggl\langle V^\text{ext}(P^k;\theta), \frac{\partial V^\text{ext}}{\partial \theta} (P^k; \theta) \biggr\rangle \\ +  \Re \langle  \lambda^1, (\nabla_{\theta} U^0) P^0 {U^0}^\dagger \rangle + \Re \langle  \lambda^1, U^0 P^0 (\nabla_{\theta} {U^0}^\dagger) \rangle \\
+ \Re \sum_{k=1}^{K-1} \langle  \lambda^{k+1}, (\nabla_{\theta} U^k) P^{k-1} {U^k}^\dagger \rangle \\
+ \Re \sum_{k=1}^{K-1} \langle  \lambda^{k+1}, U^k P^{k-1} (\nabla_{\theta} {U^k}^\dagger) \rangle.
\end{multline}
We have explained above how to compute $\nabla_{\theta} V^\text{ext}$.  The only remaining Jacobian is, for $k \geq 1$,
\begin{multline}
\label{eqn:dUdtheta}
\nabla_{\theta} U^k_{a,b} = \sum_{j,\ell} \frac{ \partial }{\partial Z_{j,\ell}} (\exp Z)_{a,b}  \biggr|_{Z = - 2 i \Delta t H(P^k, k \Delta t)} \\
\cdot (-2i \Delta t) \nabla_{\theta} H_{j,\ell}(P^k, k \Delta t; \theta)
\end{multline}
When $k=0$, we replace all instances of $2i$ with $i$, mirroring the difference between (\ref{eqn:MMUTfirststep}) and (\ref{eqn:MMUTothersteps}). For gradients of ${U^k}^\dagger$, we send $i \mapsto -i$ as explained above.  

One can check that, along feasible trajectories, the gradient (\ref{eqn:dLdtheta}) equals the gradient of the discrete objective function (\ref{eqn:discobj}) with respect to $\theta$.  Thus the basic idea is to use the gradient (\ref{eqn:dLdtheta}) as a descent direction to compute the next iterate $\theta^{(m+1)}$, repeating until convergence is achieved.  We have implemented the discrete-time formulation detailed above, using JAX \cite{jax2018github} to handle the neural network and associated automatic differentiation.  We use our implementations of the discrete objective (\ref{eqn:discobj}) and the gradient (\ref{eqn:dLdtheta}) in conjunction with a recently developed trust region quasi-Newton optimizer \cite{trophy2022}.

The quasi-Newton optimizer uses previously computed gradients to build an approximation of the Hessian; in particular, the L-SR1 (limited-memory symmetric rank-$1$) update is used \cite{trophy2022}.  This update guarantees that the resulting Hessian is symmetric but does not enforce positive definiteness.  We view this as a feature: using a neural network will almost certainly lead to an objective whose true Hessian is not positive definite.

\begin{figure*}[t]
\begin{center}
\subfloat[$\htwo$: optimal control signal $a(t)$]{\includegraphics[height=1.875in,trim={0 7pt 0 10pt},clip]{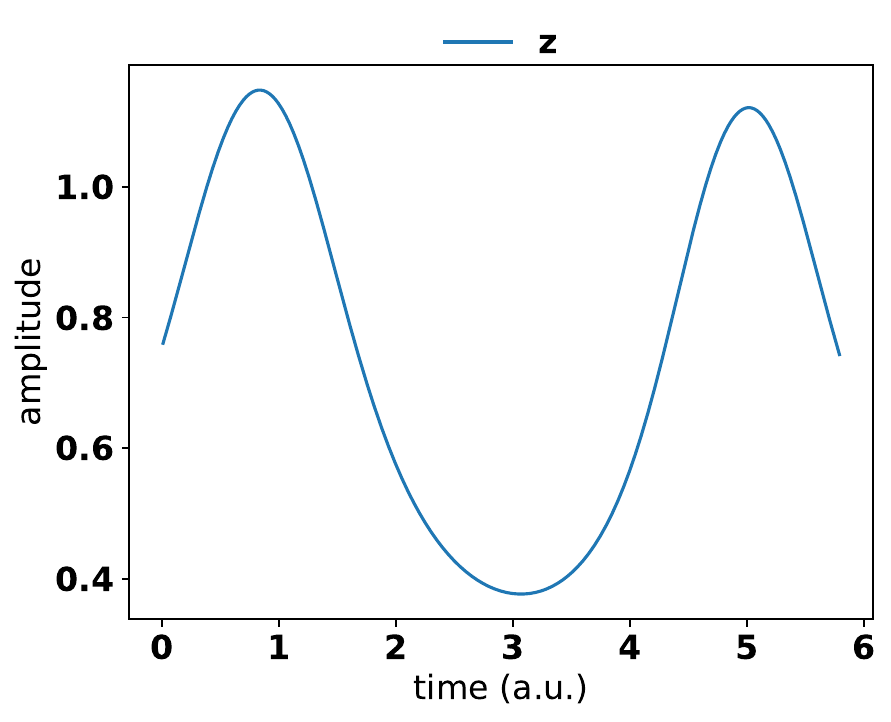}}
\subfloat[$\htwo$: controlled trajectory $P(t)$]{\includegraphics[height=1.875in,trim={0 7pt 0 10pt},clip]{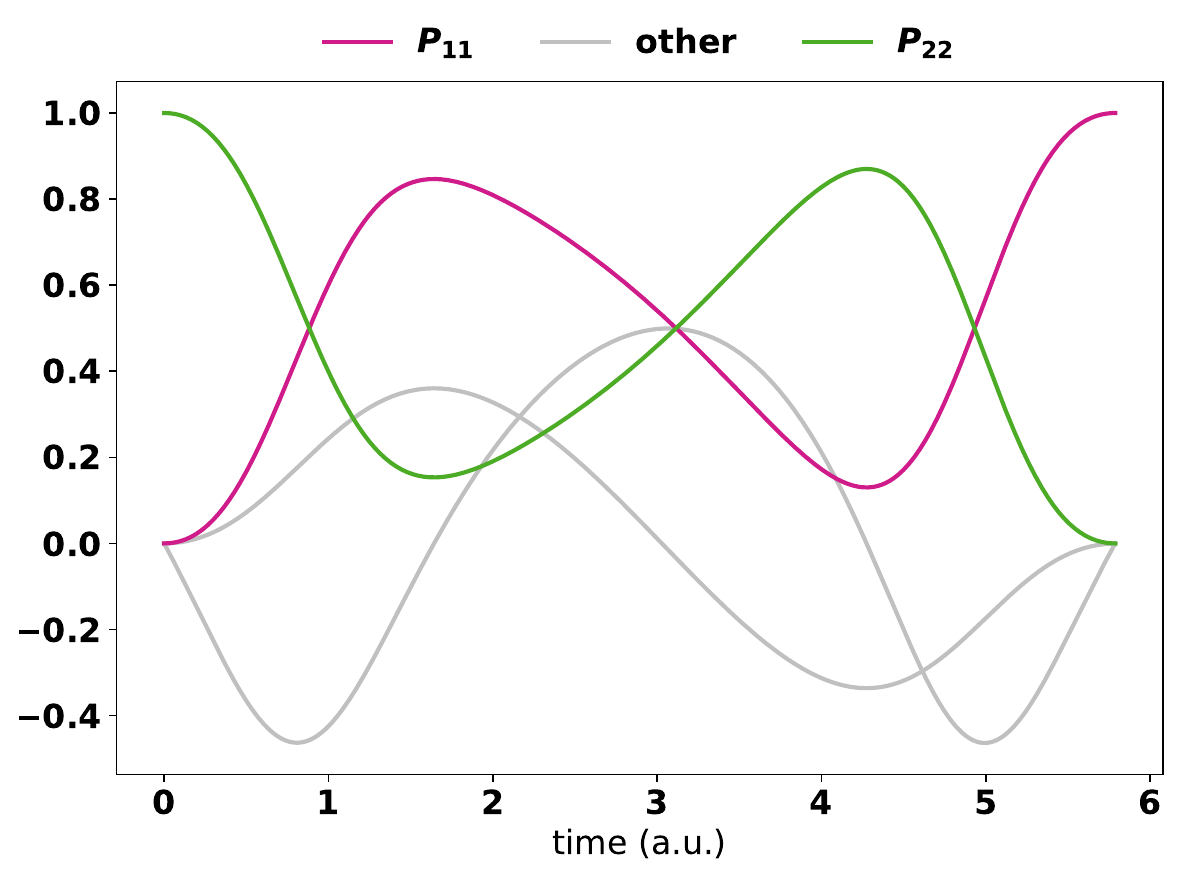}}\\ \vspace{-0.3cm}
\subfloat[$\heh$: optimal control signal $a(t)$]{\includegraphics[height=1.875in,trim={0 7pt 0 10pt},clip]{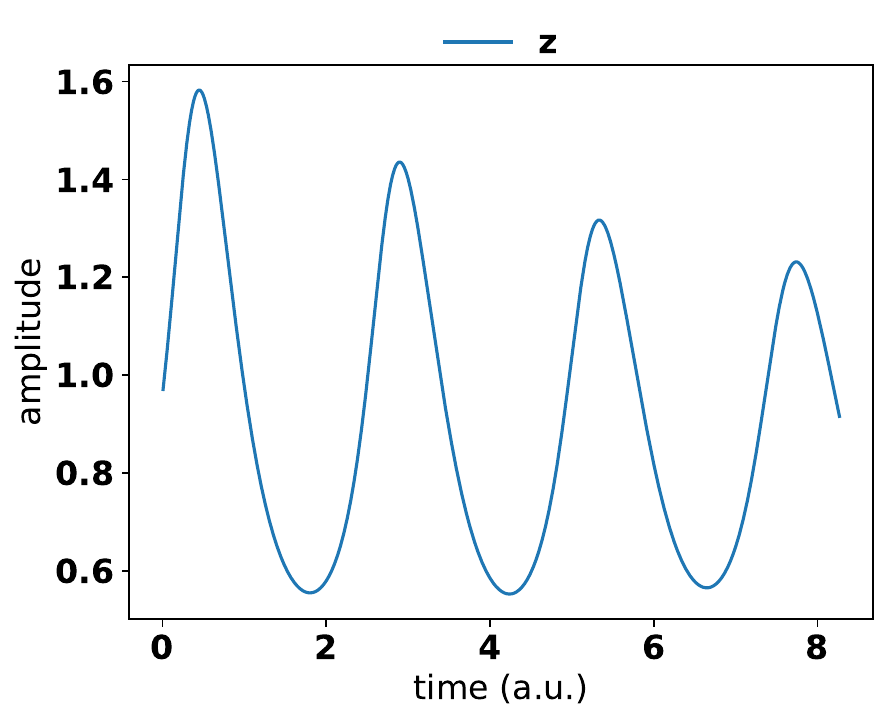}}
\subfloat[$\heh$: controlled trajectory $P(t)$]{\includegraphics[height=1.875in,trim={0 7pt 0 10pt},clip]{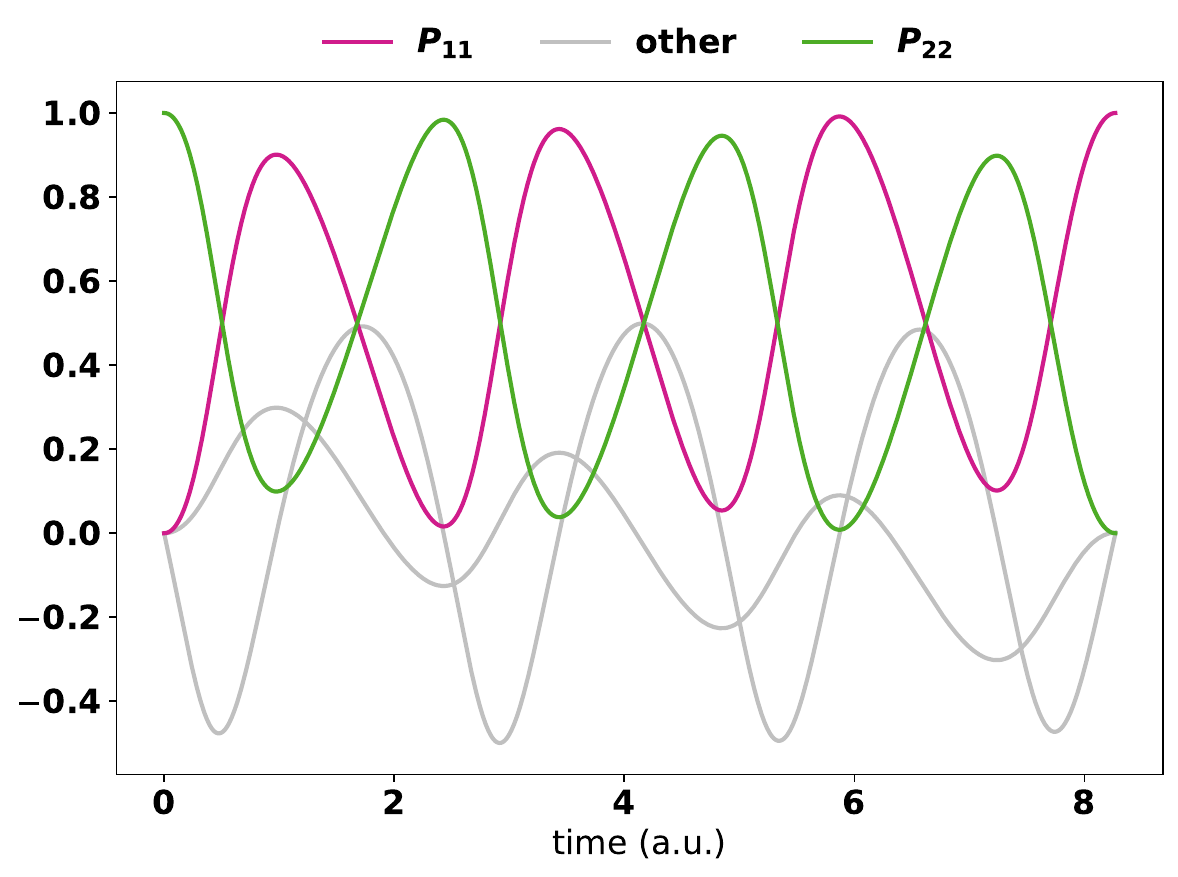}}
\end{center}
\caption{For two $2$-electron systems ($\htwo$ and $\heh$ in the STO-3G basis set), we apply the methods from this paper to learn small neural network controls $a(P)$ that drive each system (in turn) from $P_0 = \diag([0,1])$ to $P_T = \diag([1,0])$.  We then solve (\ref{eqn:tdhf}) for $P(t)$, using the Hamiltonian (\ref{eqn:hfham}) and external field (\ref{eqn:feedbackactuation}) that incorporates the control $a(P)$.  The distinct real and imaginary parts of $P(t)$ are plotted in the right panels.  The control signals $a(t) = a(P(t))$ are plotted in the left panels. Both controls achieve their desired targets with mean absolute error less than $10^{-2}$.  For these systems, the $x$ and $y$ dipole moment matrices are zero; only the control amplitude $a(t)$ in the $z$ direction matters.}
\label{fig:h2hehresults}
\end{figure*}

\section{RESULTS}
Here we give numerical results for the TDHF optimal control problem and solution methods described above.  We focus on three molecular systems: $\htwo$, $\heh$, and $\lih$ all in the STO-3G basis set.  STO-3G is a minimal basis set that uses a linear combination of three Gaussians to describe a Slater type orbital \cite[\S 3.6.2]{szabo2012modern}.  For each system, we use the Gaussian electronic structure code \cite{GaussianDV} to compute the field-free Hamiltonian $H^\text{AO}_0(P)$ in the AO basis, molecular orbital coefficients, overlap matrices, and dipole moment matrices.  Using this information, we form the field-free Hamiltonian $H_0(P)$ in the canonically orthogonalized (CO) basis \cite[\S 3.4.5]{szabo2012modern}; all remaining calculations are conducted in the CO basis.

\subsection{$\htwo$ and $\heh$}
For these $2$-electron systems, all Hamiltonians and density matrices are of size $2 \times 2$.  Symmetries dictate that only the $M_3$ dipole moment matrix (corresponding to the $z$ direction along the bond axis) is active.  Hence the control satisfies ${\mathbf{a}}(t) \in {\mathbb{R}}^1$ for all $t$. 
 For our neural network model of ${\mathbf{a}}({\mathbf{p}})$, we use a dense feedforward network with an input layer of size $N^2 = 4$, two hidden layers with $4$ units each, and an output layer of size $1$.  The total number of parameters (neural network weights and biases) in this model is $45$.  For the last layer, we use an identity activation; for prior layers, we use the softplus activation function to guarantee smoothness.

We choose $P_0 = \diag([0,1])$ and $P_T = \diag([1,0])$.  For both $2$-electron systems, all density matrices must have trace equal to $1$ at all times.  When we further constrain these matrices to be idempotent, we see that there are only two choices for diagonal matrices, the $P_0$ and $P_T$ matrices above.  This is why we have chosen them as our initial and final target states.

For both systems, we use the time step $\Delta t = 8.268 \times 10^{-3}$ a.u. (atomic units).  For $\htwo$, we use a final time of $T = 700 \Delta t$; for $\heh$, we use $T = 1000 \Delta t$.  These values were chosen after running trial trajectories with sinusoidal applied fields, to get a sense of these systems' natural time scales.

For both systems, we solve the optimal control problem repeatedly using $\rho = 10^4$ in the objective (\ref{eqn:discobj}) and random Glorot initializations \cite{pmlr-v9-glorot10a} for $\theta^{(0)}$.  \emph{For all runs, we terminate if the mean absolute difference between $P^K$ and $P_T$ is less than $10^{-2}$.}  For $\htwo$, we needed $24$ runs to obtain $10$ final results that satisfy the termination criterion within $100$ iterations; for $\heh$, we needed $14$ runs.

For $\htwo$, among the $10$ final results, there is a clear winner in terms of both mean squared control cost and fidelity: we have plotted the corresponding optimal control input $a(t)$ and controlled trajectory $P(t)$ in the top half of Figure \ref{fig:h2hehresults}.  As the control task consists of driving the system from $P_{11}(0) =0 $ and $P_{22}(0) = 1$ to $P_{11}(T) = 1$ and $P_{22}(T) = 0$, we have highlighted $P_{11}(t)$ and $P_{22}(t)$.  The remaining paths in the right panels of Figure \ref{fig:h2hehresults} are the real and imaginary parts of $P_{12}(t)$ (plotted in light gray).  By Hermitian symmetry, the real and imaginary parts of $P_{21}(t)$ need not be plotted separately.  To obtain this solution we needed only $38$ optimizer steps; the mean squared control cost is $1.28 \times 10^{-1}$ and the mean absolute difference between $P^K$ and $P_T$ is $1.98 \times 10^{-3}$.

\begin{figure*}[t]
\begin{center}
\subfloat[Optimal control signal ${\mathbf{a}}(t)$]{\includegraphics[height=2in]{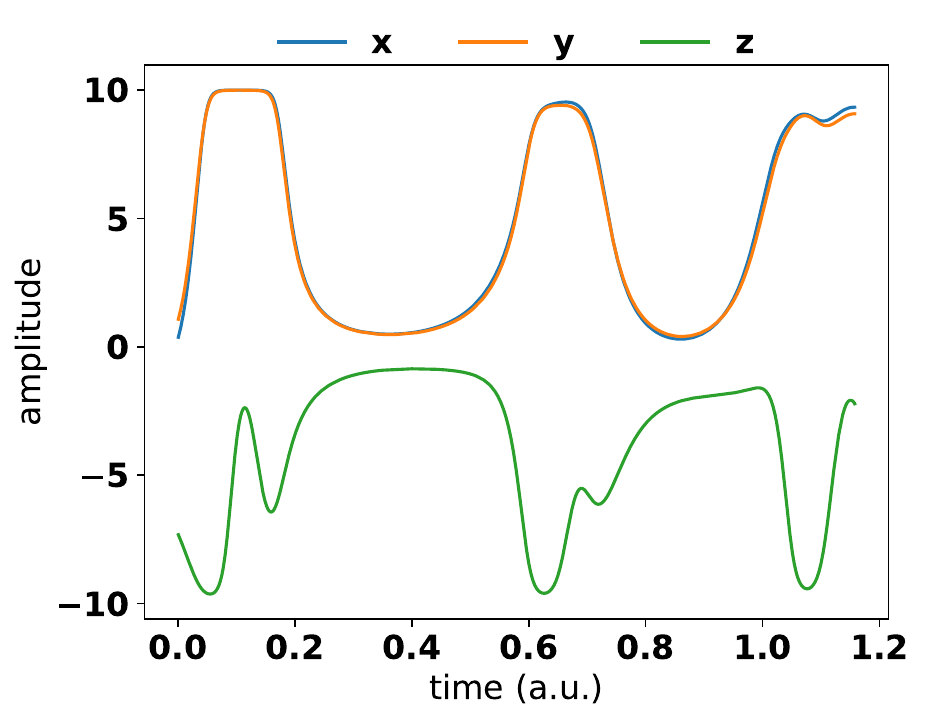}}
\subfloat[Controlled trajectory $P(t)$]{\includegraphics[height=2in]{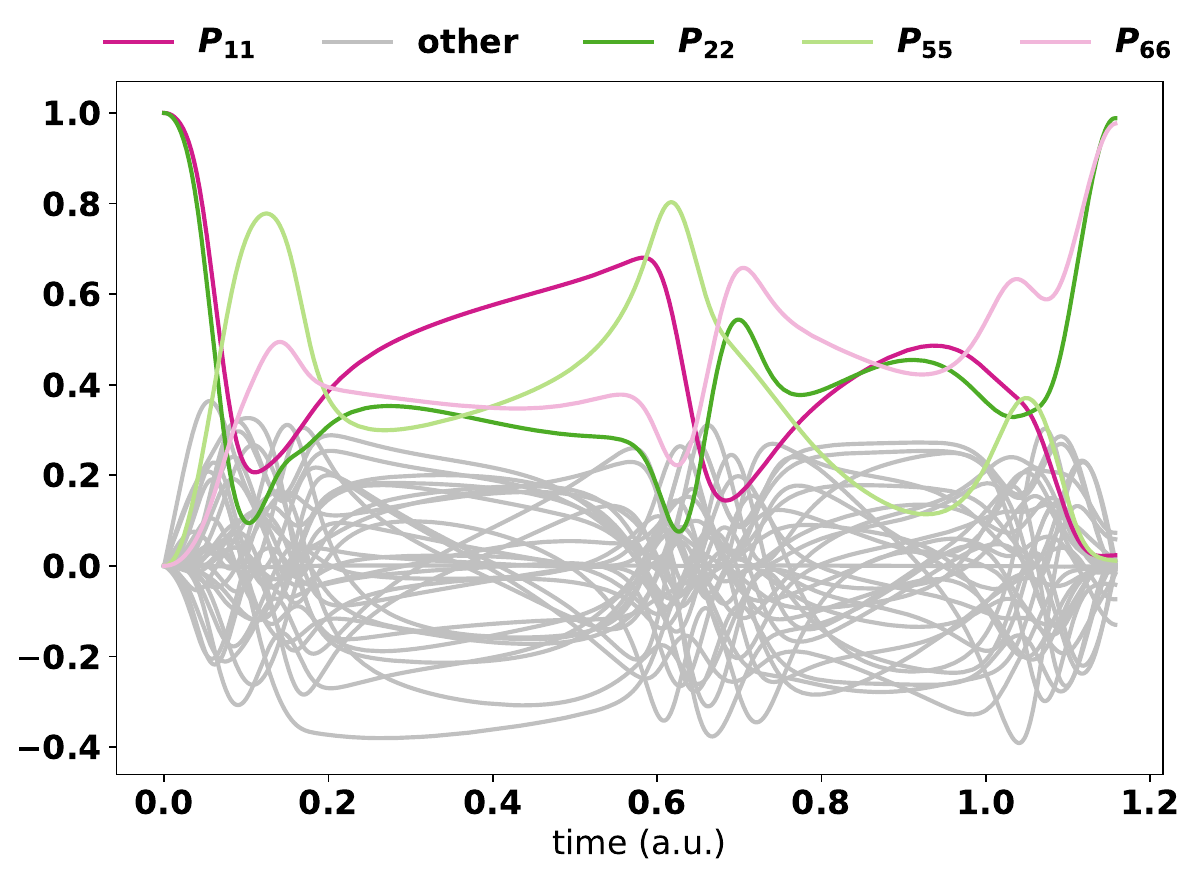}}
\end{center}
\caption{For $\lih$ in the STO-3G basis set, we apply our methods to learn a small neural network control ${\mathbf{a}}(P)$ that drives the system from $P_0 = \diag([1,1,0,0,0,0])$ to $P_T = \diag([0,1,0,0,0,1])$.  Via procedures noted in Fig. \ref{fig:h2hehresults}, we solve for $P(t)$ (right panel) and the control signal ${\mathbf{a}}(t)$ (left panel). The $x$ and $y$ components of the control signal  coincide, reflecting symmetries of the system. As the $z$ dipole moment matrix for this system does not directly couple $P_{11}$ and $P_{66}$, the control works by first moving electron density from $P_{22}$ to $P_{55}$, and then depopulating both $P_{11}$ and $P_{55}$ in favor of $P_{66}$ and $P_{22}$.  The mean absolute error between achieved and desired targets is $2.59 \times 10^{-2}$.}
\label{fig:lihresults}
\end{figure*}

For the $10$ final solutions we obtained for $\heh$, we see tradeoffs between the two terms in the objective (\ref{eqn:discobj}).  Here we describe one way to choose one solution that balances these terms.  For the $j$-th solution, let $\alpha_j$ be the mean squared control cost and $\beta_j$ be the mean absolute error between $P^K$ and $P_T$.  After rescaling all $\alpha$'s and $\beta$'s to be in $[0,1]$, we choose $j$ that minimizes $\alpha_j^2 + \beta_j^2$.  We plot in the bottom panel of Figure \ref{fig:h2hehresults} the corresponding solution, both the optimal control input and controlled trajectory.  Only $83$ optimizer steps were needed to obtain this solution; the mean squared control cost is $7.27 \times 10^{-1}$ and the mean absolute difference between $P^K$ and $P_T$ is $1.30 \times 10^{-3}$.

For both systems, at the final $\theta$ values we found through optimization, we find that $\| \nabla_{\theta} {\mathcal{J}}(\theta) \|$ is still large, on the order of $10^3$.  Analyzing this gradient further by recomputing it with $\rho=0$, we see that it primarily reflects lack of optimality with respect to the control cost.  This follows from the fact that we have driven the density matrix to a final state that is close to $P_T$, and hence we have nearly achieved optimality with respect to the fidelity.  For both systems, optimization succeeds in \emph{reducing} $\| \nabla_{\theta} {\mathcal{J}}(\theta) \|$ by at least one order of magnitude from the initial guess $\theta^{(0)}$ to the final saved $\theta$.  We hypothesize that with more iterations, the norm of the gradient will drop below desired tolerances.

\subsection{LiH}
For the $4$-electron system $\lih$, all Hamiltonians and density matrices are of size $6 \times 6$.  All three dipole moment matrices are nonzero, so $N_a = 3$ and the control satisfies ${\mathbf{a}}(t) \in {\mathbb{R}}^3$ for all $t$.  Here our neural network model of ${\mathbf{a}}({\mathbf{p}})$ is a dense feedforward network with an input layer of size $N^2 = 36$, three hidden layers with $4$ units each, and an output layer of size $3$.  The total number of parameters (neural network weights and biases) in this model is $203$.  For the final output layer, to keep the control bounded, we use an activation consisting of $z \mapsto 10 \tanh(z)$; for prior layers, to guarantee smoothness, we again use the softplus activation.

For our initial and final target states, we use two diagonal matrices:
$P_0 = \diag([1,1,0,0,0,0])$ and $P_T = \diag([0,1,0,0,0,1])$.
All density matrices have trace equal to $2$, half the number of electrons in the system.  We use $\Delta t = 8.268 \times 10^{-4}$ a.u. and final time $T = 1400 \Delta t$.

As written in the discrete objective (\ref{eqn:discobj}), the first term is $1/2$ of a sum of squares; for the $\lih$ system, we rescale this term by $1/(N^2 K)$ so that it is a \emph{mean} rather than a sum.  We use this rescaling in all subsequent gradients as well.  With this, we take $\rho = 10^3$.

We again run our optimization method repeatedly ($16$ total runs with Glorot initialization for $\theta^{(0)}$) until we obtain a set of $10$ converged solutions, using the above termination criterion and a maximum iteration limit of $1000$.  Using the same method described above to pick one solution, we arrive at the solution plotted in Figure \ref{fig:lihresults}; here the mean control cost is $13.49$ and the mean absolute difference between $P^K$ and $P_T$ is $2.59 \times 10^{-2}$.  This solution required $1000$ optimizer steps to find.  The final value of $\| \nabla_{\theta} \mathcal{J}(\theta) \|$ is large (but at least one order of magnitude smaller than at initialization), with contributions from both the control cost and the fidelity.  We hypothesize that further optimization steps would remedy this.

In the left panel of Figure \ref{fig:lihresults}, we plot the $x$, $y$, and $z$ components of the optimal amplitude ${\mathbf{a}}(t)$.  In the right panel, we have highlighted those components of the density matrix that correspond directly to the control task given by $P_0$ and $P_T$ above.  Note that the applied fields in the $x$ and the $y$ directions coincide; all applied fields are noticeably non-sinusoidal.  As compared with the smaller molecular systems described above, the controlled trajectories in Figure \ref{fig:lihresults} are noticeably more nonlinear.

\section{CONCLUSION}
This work provides proof of concept that our adjoint state method can be used to learn very small neural network controls that drive (\ref{eqn:tdhf}) from a prescribed initial state $P_0$ towards a desired target $P_T$.  This opens the door to at least three areas of future work: (i) improving the optimality of the solutions, (ii) comparison of our TDHF optimal control solutions against optimal control solutions for more exact models, and (iii) physical parametrization of the controls and interpretation of the solutions.

\section*{ACKNOWLEDGMENT}
This research was sponsored by the Office of Naval Research (Grant Number W911NF-23-1-0153) and by the US National Science Foundation (DMS-1840265).  The views and conclusions contained in this document are those of the authors and should not be interpreted as representing the official policies, either expressed or implied, of the Army Research Office or the U.S. Government. The U.S. Government is authorized to reproduce and distribute reprints for Government purposes notwithstanding any copyright notation herein.  This research used resources of the National Energy Research Scientific Computing Center (NERSC), a U.S. Department of Energy Office of Science User Facility located at Lawrence Berkeley National Laboratory, operated under Contract No. DE-AC02-05CH11231 using NERSC awards BES-m2530 and ASCR-m4577. We also acknowledge computational time on the Pinnacles cluster at UC Merced, supported by NSF OAC-2019144.

\section*{APPENDIX}
\label{sect:appendix}
\noindent For the Frobenius inner product (\ref{eqn:froip}), we can establish
\begin{align}
\langle A, B \rangle &= \sum_{i=1}^n \sum_{j=1}^n (A^\dagger)_{ij} B_{ji} = \sum_{i,j} \overline{A_{ji}} B_{ji}, \\
\label{eqn:conj}
\overline{ \langle A, B \rangle } &= \sum_{i,j} A_{ji} \overline{B_{ji}} = \langle \overline{A}, \overline{B} \rangle, \text{ and } \\
\label{eqn:fro}
\langle A, A \rangle &= \| A \|_F^2 = \sum_{i,j} | A_{ji} |^2,
\end{align}
the squared Frobenius norm of $A$.  From (\ref{eqn:conj}), we have
\begin{equation}
\label{eqn:conj2}
\Re \langle A, B \rangle = \Re \overline{ \langle A, B \rangle } = \Re \langle \overline{A}, \overline{B} \rangle.
\end{equation}
Using $\tr P Q = \tr Q P$, we can show
\begin{align}
\langle A, B \rangle &= \tr A^\dagger B = \tr B A^\dagger = \langle B^\dagger, A^\dagger \rangle, \\
\label{eqn:mat12}
\langle A, Q B \rangle &= \langle Q^\dagger A, B \rangle \ \text{ and } \ \langle A, B Q \rangle = \langle A Q^\dagger, B \rangle.
\end{align}
Suppose that $\xi$ is a $4$-index tensor.  We use $:$ to denote the double-index contractions
\begin{equation}
\label{eqn:contractions}
(\xi : A)_{ij} \! = \! \sum_{k,\ell} \xi_{ijk\ell} A_{k\ell}, \ (A : \xi)_{ij} \! = \! \sum_{k,\ell} A_{k\ell} \xi_{k\ell ij}.
\end{equation}
The result of contracting $\xi$ against a matrix is a matrix.  We can then see that
\begin{multline}
\label{eqn:tens}
\langle A, \xi : B \rangle = \sum_{i,j} \overline{A_{ji}} (\xi : B)_{ji} = \sum_{\ell,k} \sum_{j,i} \overline{ A_{ji} \overline{\xi_{jik\ell}} } B_{k \ell} \\
= \sum_{\ell,k} \overline{(A : \overline{\xi})_{k \ell}} B_{k \ell} = \langle A : \overline{\xi}, B \rangle
\end{multline}


{\small %
}

\end{document}